\documentclass[11pt,reqno]{amsart}
\usepackage{amsmath,amssymb,graphicx,amscd,hyperref, amsfonts,psfrag,fancybox,qsymbols,indentfirst}
\usepackage{mathrsfs}
\usepackage{dsfont}
\usepackage[latin1]{inputenc}
\newcommand{\K}{{\kappa}}

 \DeclareMathOperator{\E}{\mathbb{E}}      
 \DeclareMathOperator{\Var}{\mathrm{Var}}      
 \DeclareMathOperator{\cov}{cov}             
\DeclareMathOperator{\Tr}{Tr}
\DeclareMathOperator{\tr}{tr}

\DeclareMathOperator{\M}{ M\ddot{o}b}

\DeclareMathOperator{\W}{Wg}
\DeclareMathOperator{\Wl}{W\Lambda}

\newcommand{\convlaw}{\overset{\mbox{\rm \scriptsize law}}{\longrightarrow}}

\newtheorem{thm}{Theorem}[section]
\newtheorem{cor}[thm]{Corollary}
\newtheorem{lem}[thm]{Lemma}
\newtheorem{prop}[thm]{Proposition}
\newtheorem{rem}[thm]{Remark}
\theoremstyle{definition}
\newtheorem{defn}[thm]{Definition}
\theoremstyle{remark}

\def \be{\begin{eqnarray*}}
\def \ee{\end{eqnarray*}}
\def \ben{\begin{eqnarray}}
\def \een{\end{eqnarray}}

\def\QED{\hfill\vrule height 1.5ex width 1.4ex depth -.1ex \vskip20pt}

\numberwithin{equation}{section}

\begin{document}
\title[Truncations, traces and bivariate Brownian bridge]
{Truncations of Haar distributed matrices, traces  
and bivariate Brownian bridge}
\date{\today}

\author{C. Donati-Martin}
\address{Université Paris 6 and CNRS, UMR 7599, Laboratoire de Probabilités et Modèles Aléatoires, Case courrier 188
4 place Jussieu, F-75252 Paris Cedex 05}
\email{catherine.donati@upmc.fr}

\author{A. Rouault}
 \address{Universit\'e  Versailles-Saint Quentin, LMV UMR 8100,
 B\^atiment Fermat, 45 avenue des Etats-Unis,
F-78035 Versailles Cedex}
 \email{alain.rouault@math.uvsq.fr}


\begin{abstract}
Let $U$ be a Haar distributed 
 matrix in $\mathbb U(n)$ or $\mathbb O(n)$. We show that after centering the two-parameter process 
\[W^{(n)} (s,t) = \sum_{i \leq \lfloor ns \rfloor, j \leq \lfloor nt\rfloor} |U_{ij}|^2\]
converges in distribution to the bivariate tied-down Brownian bridge. 
\end{abstract}


\maketitle
\section{Introduction}
Let $\sigma$ be a random permutation uniformly distributed on the symmetric group ${\mathcal S}_n$.
Define for $p, q \leq n$
$$ X_{p,q}^{(n)} = 
\# \{1\leq i \leq p \ ,\ \sigma(i) \leq q \}.$$
In \cite{chapuy}, G. Chapuy proved that a suitable normalization of  $X_{p,q}^{(n)} $ converges in distribution to the bivariate tied-down Brownian bridge. 
Note that $X_{p,q}^{(n)}  = \Tr(\Sigma_{p,q} \Sigma_{p,q}^\star)$ where $\Sigma_{p,q}$ is the truncated matrix of size $p\times q$ of $\Sigma$, the permutation matrix associated with $\sigma$, and $^\star$ means adjoint.
In this paper, we prove a similar  result when the symmetric group is replaced by the unitary group or the orthogonal group, equipped with the Haar measure. 

Let $U$ be a Haar unitary, resp. orthogonal,  matrix in $\mathbb U (n)$, resp. in $\mathbb O (n)$. We consider,   for $p \leq n$ and $q\leq n$, 
the upper-left $p\times q$ submatrix $V_{p,q}$ and the $p\times p$
Hermitian matrix 
\[H_{p,q} = V_{p,q}V^\star_{p,q}\,.\] We are interested in the asymptotic behavior of
\begin{equation}
\label{tpq}T_{p, q} = \Tr H_{p,q} = \sum_{i \leq p, j\leq q} |U_{i,j}|^2\,.\end{equation} 
Setting
\[Y^{(n)}_{p,q} = T_{p, q} - \mathbb E T_{p, q}\,,\]
we define a sequence of two-parameter processes $W^{(n)}$  by
\[W^{(n)} := \left( Y^{(n)}_{\lfloor ns\rfloor,\lfloor nt\rfloor} , s, t \in [0,1]\right)\]

Chapuy used the space $C([0,1]^2)$ completing its process 
 in such a way that it is continuous and affine on each closed "lattice triangle". We prefer using the 
 multidimensional generalization of Skorokhod space
  $D([0,1]^2)$ given by \cite{bickel1971convergence}.  It consists of functions
from $[0, 1]^2$ to $\mathbb R$ which are at each point right continuous (with respect
to the natural partial order of $[0, 1]^2$) and admit limits in all "orthants". 
The space  $D([0,1]^2)$ is endowed with the topology of Skorokhod (see \cite{bickel1971convergence} for the definition).

\noindent
Our main result is the following
\begin{thm} \label{maintheo}
The process $W^{(n)}$  converges in distribution in $D([0,1]^2)$ to a tied-down  Brownian bridge $\sqrt{ \frac{2}{\beta}} W^{(\infty)}$
where $W^{(\infty)}$  is a centered continuous Gaussian process on $[0,1]^2$ of covariance
$$\E[ W^{(\infty)}(s,t) W^{(\infty)}(s',t')] = (s \wedge s' -s s')(t\wedge t' - tt'),$$
 $\beta = 2$ in the unitary case and $\beta = 1$ in the orthogonal case.
\end{thm}

Previous works are related to our problem. First, Borel in 1906 shows that for a uniformly distributed point on the $(n-1)$-dimensional (real) sphere, the scaled first coordinate converges in distribution to the standard normal. Since that time, many authors studied the entries and partial traces of matrices from the orthogonal and  unitary group. In particular Diaconis and d'Aristotile (\cite{Aristo1, Aristo2})  
 proved 
  that the sequence of one-parameter processes
\[\left\{\sum_{i=1}^{\lfloor ns\rfloor}  U_{ii} \ , \ s \in [0,1]\right\}_n\] converges in distribution to the complex Brownian motion.
Besides,   
  Silverstein \cite{Silver1}
proved  that for $q$ fixed, the sequence of one-parameter processes \begin{equation}
\label{siver}\left\{n^{1/2}(\sum_{i=1}^{\lfloor ns\rfloor} |U_{iq}|^2 - s) \ , \ s\in [0,1]\right\}_n\end{equation}  
converges in distribution to the Brownian bridge. 

In the paper \cite{Silver2}, Silverstein  discussed the similarity between the matrix of eigenvectors of a (real) sample covariance matrix and a Haar distributed orthogonal matrix, with a one-parameter parameter process analogous to (\ref{siver}). We also refer to \cite[Chapter 10]{Baibook} for the behavior of eigenvectors of sample covariance matrices and universality conjectures. 
To extend this study and after reading a first draft of our result, Djalil Chafai conjectured that if  $M$ is a $n\times n$ matrix with i.i.d. entries  
 having the same four first moments as the complex Gaussian standard 
 and if $\mathcal U$ denotes the matrix of  eigenvectors of  $N=MM^*$, then the sequence  ${\mathcal W}^{(n)}$ obtained by changing $U_{ij}$ into 
${\mathcal U}_{ij}$ converges to the tied-down Brownian bridge as in Theorem \ref{maintheo}. 
At the moment of posting  
the present version,  we are aware that Florent Benaych-Georges (\cite{FBG}) answered positively to the conjecture when $N$ is a Wigner matrix, under a fourth moment hypothesis.

The rest of the paper is organized as follows. In section 2, we introduce the basic notions on permutations, partitions, classical cumulants. The paper is then split into two parts. Section 3 deals with the unitary case and Section 4 is devoted to the orthogonal case. In each part, we develop the combinatorics associated with the group\footnote{We will deal with the symplectic case in a forthcoming paper.} 
(the Weingarten function and the associated cumulants).

In particular,  we state a formula for the cumulants of variables of the form $X = \Tr( AUBU^\star)$ for deterministic matrices $A,B$ of size $n$. In the unitary case, the formula follows from the results of \cite{mingo2007second}. 
We then apply the above formula to the computation of  the second and fourth cumulants of $T_{p,q}$. The proof of Theorem \ref{maintheo} is then divided in two parts: tightness of the family of distributions of $W^{(n)}$ and convergence of the finite dimensional laws. To prove the tightness, we use a criterion of Bickel and Wichura for two-parameter processes, with the help of the estimates  obtained in section 3.2, resp 4.3. The convergence of the finite dimensional distributions to  Gaussian distributions relies on the computations of their cumulants and their asymptotics. The expression of the cumulants follows from the previous section and their limit follows from the  asymptotics of cumulants of  unitary, resp. orthogonal, Weingarten functions, obtained in \cite{Co}, resp. \cite{CoSn}.
Let us mention that simultaneously and independently of the present paper, precise computations in the orthogonal case are presented in
\cite{Redelmeier1} and \cite{Redelmeier2}. 
In section 5, we give complementary remarks and connections with other problems.
\section{Preliminaries}
For $n$ a positive integer we set $[n] := \{1,2, \cdots, n\}$. 
\subsection{Partitions, permutations.}

We call $A= \{A_1, \cdots, A_s\}$ a partition of $[k]$ if the $A_i \ (1\leq i \leq s)$ are pairly disjoint, non-void subsets of $[k]$ such that $A_1 \cup\cdots\cup A_s = [k]$. We call $A_1, \cdots, A_s$ the blocks of $A$. The number of blocks of $A$ is denoted by $\#(A)$. The set of all partitions of $[k]$ is denoted by ${\mathcal P}(k)$. One partition $A$ is said to refine another $B$, denoted $A \leq B$ provided every block of $A$ is contained in some block of $B$. Given two partitions $A$ and $B$, $A\wedge B$ (resp. $A\vee B$) is the largest (resp. smallest) partition which refines (resp. is refined by) both $A$ and $B$. Under these operations, the partially ordered set ${\mathcal P}(k)$ is a lattice. We denote by $1_k$   the largest partition of $[k]$ (one-block partition), and by $0_k$ the smallest one ($k$-blocks partition). For $A, B \in {\mathcal P}(k)$ with $A \leq B$ we denote by $[A,B]$ the interval 
\[[A,B] = \{  C \in {\mathcal P}(k)\ |\ A \leq C\leq B\}\,.\]
Since we will make some use of the Möbius function for partitions, let us just recall (\cite{Stanley} Sect 3.6) that for two real  functions 
$f,g$ defined on $\{(A,B) \in  {\mathcal P}(k)\times {\mathcal P}(k); A \leq B \}$, 
we have: 
\begin{equation}
\label{2.1.1}f(A,B) = \sum_{C\in [A,B]} g(A,C) 
\end{equation}
if and only if
\begin{equation}
\label{2.1.2}g(A,B) = \sum_{C\in [A,B]} \M(C,B)f(A, C) \ \ \ \ 
\end{equation}
with 
\[\M(C,B) := \prod_i \left((-1)^{i-1}(i-1)!\right)^{p_i}\] 
where  $p_i$ is the number of blocks of $B$ that contain exactly $i$ blocks of $C$.

Let ${\mathcal S}_k$ be the set of permutations on $k$ elements. 
With $\sigma \in {\mathcal S}_k$, we associate the set ${\mathcal C}(\sigma)$ of its cycles, whose number is denoted by $\#(\sigma)$. 
  We denote by $0_\sigma$ the partition whose blocks are the cycles of $\sigma$, or when the context is clear, just by $\sigma$.
 For $\pi \in {\mathcal S}_{k}$, a partition $A=(A_1, \ldots A_l)$ of $[k]$ is called $\pi$-invariant  if $\pi$ leaves invariant each block $A_i$ that is $0_\pi \leq A$ (which we just write $\pi \leq A$). \\

Finally,  we define ${\mathcal M}_{2k}$ as  the set of pairings of $[2k]$, i.e. of partitions where each block consists of exactly two elements. 
 It is then convenient to encode  the set $[2k]$ by
 \[[2k] \cong \{1, \ldots, k, \bar 1, \ldots, \bar k\}\,.\]
   Given two pairings $p_1, p_2$, we define the graph $\Gamma(p_1, p_2)$ as follows. The vertex set is $[2k]$ and the edge set consists of the pairs of $p_1$ and $p_2$.
   Let ${\rm loop}(p_1, p_2)$ the  number  of connected components of $\Gamma(p_1,p_2)$. 

 \medskip
 
  
\subsection{Cumulants}
 For $r \geq 1$, $\K_r$ denotes the classical cumulant of order $r$ (see \cite{Sh},  \cite{mingo2007second} p.215). It is 
 a  multilinear function of $r$ variables defined as follows: if $a_1, \ldots, a_r$ are random variables,
$$\K_r(a_1, \ldots ,a_r) = \sum_{C \in {\mathcal P}(r)} \M(C, 1_r) \E_C(a_1, \ldots a_r)$$
where for $C = \{C_1, \ldots , C_k\} \in {\mathcal P}(r)$,
\begin{equation}
\label{defEC}  \E_C(a_1, \ldots , a_r) = \prod_{i=1}^r \E(\prod_{j\in C_i} a_j)\,.
\end{equation}
More generally, relative cumulants are defined, for $A \leq B\in {\mathcal P}(r)$ as
\begin{equation}
\label{relcum}\kappa_{A, B} (a_1, \cdots, a_r) = \sum_{C\in [A,B]} \mathbb E_C (a_1, \cdots, a_r) \M(C, B)\,.\end{equation}
From the equivalence between (\ref{2.1.1}) and (\ref{2.1.2}) we have, for any $A \leq B \in {\mathcal P}(r)$
\begin{equation}
\label{relcumrev}\mathbb E_{B} (a_1, \cdots, a_r) = \sum_{C\in [A,B]} \kappa_{A, C} (a_1, \cdots, a_r)\,.
\end{equation}

\subsection{Matrices}
For a matrix $M= (M_{ij})_{i,j \leq n}$, we denote by $\Tr$ the trace and by $\tr$ the normalized trace
$$\Tr(M) = \sum_{i=1}^n M_{ii}, \qquad  \tr(M) = \frac{1}{n} \sum_{i=1}^n M_{ii}.$$
 For $\pi \in {\mathcal S}_k$ and $M = (M_1, \ldots M_k)$  a k-tuple of $n\times n$ matrices, we set
$$ \Tr_\pi(M) = \Tr_\pi(M_1, \ldots , M_k) = \prod_{C \in {\mathcal C}(\pi) } \Tr(\prod_{j \in C} M_j).$$

  Let $s$ be a fixed integer and let $\{M_1, \ldots ,M_s\}_n$ be a sequence of $n\times n$ deterministic matrices. We say that 
$\{M_1, \ldots M_s\}_n$ has a limit distribution if there exists  a non commutative probability space $({\mathcal A}, \varphi)$ and $a_1, \ldots a_s \in {\mathcal A}$ such that for any polynomial $p$ in $s$ non commuting variables,
$$ \lim_{n \rightarrow \infty} \tr(p(M_1, \ldots, M_s))  = \varphi(p(a_1, \ldots, a_s)).$$

\section{The unitary group}
\subsection{Preliminary remarks: Some moments}
\label{prelimunit}
Let $U$ be a Haar distributed matrix on $\mathbb U(n)$ the unitary group of size $n$.
We have the following important relations (see \cite{HiaiP}, Proposition 4.2.3 and \cite{Ji}). If $U_{i,j}$ 
is the generic element of $U$, 
 then the random variable $|U_{i,j}|^2$
follows the  
 beta distribution on $[0,1]$ with  parameter $(1, n-1)$ of density $(n-1)(1-x)^{n-2}$.
%
 Thus,
\ben
\label{momentX}
\mathbb E |U_{i,j}|^2 = \frac{1}{n} \ , \ \mathbb E |U_{i,j}|^4 = \frac{2}{n(n+1)} \ , \ \Var |U_{i,j}|^2 =  \frac{n-1}{n^2(n+1)}\,, 
\een
and more generally
\ben
\mathbb E |U_{ij}|^{2k} = \frac{(n-1)! k!}{(n-1+k)!}\,.
\een
If $X = |U_{i,j}|^2$ and $Y = |U_{i,k}|^2$ with $k\not = j$, then $(X,Y)$ follows the Dirichlet distribution on $\{0\leq x, y , x+y \leq 1\}$ with parameters  
$(1, 1, n-2)$ of density $(n-1)(n-2)(1-x-y)^{n-3}$. Thus
\ben \label{2.2}
\mathbb E \left(|U_{i,j}|^2|U_{i,k}|^2\right) = \frac{1}{n(n+1)}\,.
\een
Besides, if $i\not= k, j \not= \ell$,
\ben \label{2.3}
\mathbb E \left(|U_{i,j}|^2|U_{k,\ell}|^2\right) = \frac{1}{n^2 - 1}\,. 
\een
From these relations, we can compute the first moments of $T_{p,q}$, defined in (\ref{tpq}). 
\begin{prop} \label{propvar}
The mean and the variance of $T_{p,q}$ are given by:
\ben
\label{mean}
\mathbb E T_{p,q} = \sum_{i\leq p, j \leq q} \mathbb E |U_{ij}|^2 = pq \mathbb E |U_{11}|^2 = \frac{pq}{n}\,. 
\een
and 
\ben \label{var0}
\Var T_{p,q} = pq \frac{n^2 - n(p+q) + pq}{n^2(n^2-1)}. \een
Assume that $p/n \rightarrow s$, $q/n \rightarrow t$, then,
$$ \lim_n \frac{1}{n} \mathbb E T_{p,q} = st \ , \ \lim_n \Var T_{p,q}  = st(1-(s+t) +st) = st(1-s)(1-t).$$
\end{prop}
{\bf Proof:} \begin{eqnarray*}
\mathbb E T^2_{p,q} &=&\sum_{i,k\leq p, j,l \leq q} \mathbb E |U_{ij}|^2 |U_{kl}|^2 \\
&=& \sum_{i\leq p, j \leq q} \mathbb E |U_{ij}|^4 + \sum_{i\leq p, j\not= l \leq q } \mathbb E |U_{ij}|^2 |U_{il}|^2  \\
&& \qquad
+ \sum_{i\not= k \leq p, j \leq q} \mathbb E |U_{ij}|^2 |U_{kj}|^2 + \sum_{i\not=k \leq p,j \not= l \leq q } \mathbb E |U_{ij}|^2 |U_{kl}|^2 \\
&=& pq \frac{2}{n(n+1)} + pq(q-1)  \frac{1}{n(n+1)} \\
&& \qquad + p(p-1)q  \frac{1}{n(n+1)} + p(p-1)q(q-1)  \frac{1}{n^2-1} \\
&=& pq \left( \frac{p+q}{n(n+1)} +  \frac{(p-1)(q-1)}{n^2-1}\right)\,.
\end{eqnarray*}
This yields \eqref{var0}.

\begin{rem}
An easy consequence of the above Proposition is
\[\lim_n \frac{1}{n} T_{\lfloor ns \rfloor, \lfloor nt\rfloor} = st\]
in probability. Actually the convergence is uniform in $s, t \in [0,1]$ (see section 5).
\end{rem}

\subsection{Combinatorics for the unitary group}
Let ${\mathbb U}(n)$ denote the unitary group of size $n$ endowed with the Haar probability measure. The generic element will be denoted by $U$ and its $(i,j)$ coefficient by $U_{ij}$.
In \cite{CoSn}, Collins and Sniady proved the following integration formula on ${\mathbb U}(n)$, see also \cite{Co}. 
Let ${\mathcal M}_{2k}^U$ denote the set of pairings of $[2k]$, pairing each element of $[k]$ with an element of $[\bar k]$. Let $G^{{\mathbb U}(n)}$ be the Gram matrix\footnote{The term of Gram matrix comes from the theory of  representations of groups and algebras, see \cite{CoMa}.}
$$ G^{{\mathbb U}(n)} = (G^{{\mathbb U}(n)}(p_1, p_2))_{p_1,p_2 \in {\mathcal M}_{2k}^U} := (n^{{\rm loop}(p_1, p_2)})_{p_1, p_2\in {\mathcal M}_{2k}^U}\,.$$
Then the unitary Weingarten matrix 
 $\W^{{\mathbb U}(n)}$
 is the pseudo inverse of $G^{{\mathbb U}(n)}$.

\begin{prop} \label{propWeing-unit}
For every choice of indices ${\mathbf i} = (i_1, \ldots, i_{k}, i_{\bar 1}, \ldots, i_{\bar k})$ and 
 ${\mathbf j} = (j_1, \ldots, j_{k}, j_{\bar 1}, \ldots, j_{\bar k})$,
\begin{equation} \label{Weing-unit1}
\mathbb E \left(U_{i_1 j_1} \ldots U_{i_k j_k} \bar U_{i_{\bar 1}j_{\bar 1}}\ldots \bar U_{i_{\bar k}, j_{\bar k}}\right) = 
\sum_{p_1, p_2 \in {\mathcal M}_{2k}^U} \delta_{\mathbf i}^{p_1}  \delta_{\mathbf j}^{p_2} \W^{{\mathbb U}(n)}(p_1, p_2)
\end{equation}
where 
 $\delta_{\mathbf i}^{p_1}$ (resp. $\delta_{\mathbf j}^{p_2})$ is equal to $1$  or $0$ if $\mathbf i$ (resp. $\mathbf j$) is constant on each pair of $p_1$ (resp. $p_2$) or not.
\end{prop}
It is clear that with each pairing $p \in {\mathcal M}_{2k}^U$ we can associate a unique $\sigma \in {\mathcal S}_k$ such that $p = \prod_{i=1}^k (i, \overline{\sigma(i)})$. It is known that if $p_1$ is associated with $\alpha$ and $p_2$ with $\beta$, then  $\W^{{\mathbb U}(n)}(p_1, p_2)$ is a function of $\beta\alpha^{-1}$ denoted by $\W(n, \beta\alpha^{-1})$, so that  (\ref{Weing-unit1}) becomes
\begin{equation} \label{Weing-unit2}
\mathbb E\left(  U_{i_1j_1} \ldots U_{i_k, j_k} \bar U_{i_{\bar 1}j_{\bar 1}}\ldots \bar U_{i_{\bar k}, j_{\bar k}}\right) 
= \sum_{\alpha, \beta \in {\mathcal S}_{k}} \tilde\delta_{\mathbf i}^{\alpha}  \tilde\delta_{\mathbf j}^{\beta} \W(n , \beta\alpha^{-1})
\end{equation}
where $\tilde\delta_{\mathbf i}^{\alpha} = 1$ if $i(s)= i(\overline{\alpha(s)})$ for every $s \leq k$ and $0$ otherwise. In particular, if $\pi \in {\mathcal S}_k$, we have
\begin{equation}
\label{Weinu}
\W(n, \pi) = \E(U_{11} \ldots U_{kk} \bar{U}_{1\pi(1)}\ldots \bar{U}_{k\pi(k)})\,.
\end{equation}
 \\
The Weingarten functions for $k=1,2$ are given by (see \cite{Co}):
\ben \label{Wg}
\W(n, (1)) &=  & \frac{1}{n} \nonumber\\ 
 \W(n, (1)(2)) &= & \frac{1}{n^2-1}\ , \ \W (n, (12)) =   - \frac{1}{n(n^2-1)}\,.
 \een
 From these equations, we can recover \eqref{2.2}, \eqref{2.3}.

 \noindent
We can now state a proposition which is a particular case  of \cite[Theorem 3.10]{mingo2007second}.

\begin{prop} \label{mainprop2}
Let $U$ be Haar distributed on ${\mathbb U}(n)$.
Let  $D= (D_1, \ldots D_k)$ and $\bar{D}= (D_{\bar{1}}, \ldots 
D_{\bar{k}})$ be two families  of deterministic matrices of size $n$.
We set, for $ 1 \leq i \leq r$,
$$X_i = \Tr(D_iUD_{\bar{i}}U^{\star})\,.$$
Then,
\begin{equation} \label{main2u}
\kappa_r(X_1, \ldots, X_r)= 
\sum_{\alpha, \beta \in {\mathcal S}_{r} } \sum_{A}  C_{\beta\alpha^{-1},  A} \Tr_{\alpha}(\bar D) \Tr_{\beta^{-1}}(D)
\end{equation}
where in the second sum $A \in {\mathcal P}(r)$ is such that
\begin{equation}
\label{condu}
\beta\alpha^{-1} \leq A \ \hbox{and}\  
  A\vee \beta \vee \alpha= 1_{r}\,,
  \end{equation}
and $C_{\sigma, A}$ are the relative cumulants of the unitary Weingarten function (see \cite{CoSn})\,.
Moreover, if the sequence $\{D, \bar D\}_n$ has a limit distribution, then for $r \geq 3$,
$$\lim_{n \rightarrow \infty} \K_r(X_1, \ldots, X_r) = 0.$$
\end{prop}
\begin{rem} When writing $\Tr_{\alpha}(\bar D)$ in \eqref{main2u}, we consider $\alpha$ as a permutation acting on $[\bar{k}]$. The formula \eqref{main2u} is not given exactly on this form in \cite{mingo2007second} but in our particular case where $X_i = \Tr(D_iUD_{\bar{i}}U^{\star})$ with the $D_i$ deterministic, the formulas are equivalent. We shall prove a similar formula in the orthogonal case.
\end{rem}

\noindent For the sake of completeness, let us give the meaning of $C_{\pi , A}$. Let $\pi \in {\mathcal S}_r$ and $a_i = U_{ii}\bar U_{i \pi(i)}$. We denote for a $\pi$ invariant partition $\Pi \in {\mathcal P}(r)$
\begin{equation}
E_\Pi (\pi) := \mathbb E_\Pi (a_1, \ldots, a_r)\,,
\end{equation} 
i.e., owing to (\ref{Weinu})
\begin{equation}
E_\Pi (\pi) = \prod_{k=1}^s \W(\pi|_{V_k})\,,
\end{equation} 
where $\Pi = \{V_1, \ldots, V_s\}$. Then, if we set
\[C_{\pi, A} := \kappa_{\pi, A}(a_1, \ldots, a_r)\,,\]
the summation formula (\ref{relcum}) yields
\begin{equation}
\label{relcumU}
C_{\pi, A} = \sum_{\pi \leq C \leq A} \W(\pi|_{V_1})\cdots \W(\pi|_{V_1}) \M(C,A)
\end{equation}
and the reverse one (\ref{relcumrev})
\begin{equation}
\label{relcumrevU}
E_C (\pi) = \sum_{A \in [\pi, C]} C_{\pi, A}\,. 
\end{equation}
We will use these formulas in the orthogonal case.
\subsection{Computations of the second and fourth cumulants of $T_{p,q}$}
\subsubsection{The covariance of $T_{p,q}$}
The fundamental remark is that
\[H_{p,q} = D_1UD_{\bar 1}U^\star\]
with $D_1 = I_p, D_{\bar 1} = I_q$,
 where $I_k$ is the matrix of projection on the $k$ first coordinates. Note that $D_1$ and $D_{\bar 1}$ are commuting projectors and that if $p/n \rightarrow s, q/n \rightarrow t$, $\{D_1, D_{\bar 1}\}_n$ has  a limit distribution with $a_1, a_2$ commuting projectors on $({\mathcal A} , \varphi)$ such that  $a_1a_2 =a_1$ if $s<t$ and $=a_2$ if $t<s$, and $\varphi(a_1) =s, \varphi(a_2) =t$. \\
Let $p,p',q, q' \leq n$. We now give an application of  Proposition \ref{mainprop2} to the computation of $\cov(T_{p,q}, T_{p',q'})= \K_2(T_{p,q}, T_{p',q'})$. This can also be done, using the computations of Section \ref{prelimunit}.\\
We set $D_{2}= I_{p'}, D_{\bar 2} = I_{q'}$ and apply formula \eqref{main2u} to $X_1= T_{p,q}$, $X_2 = T_{p',q'}$, $r=2$. 
The different possibilities for $\alpha, \beta$ and $A$ satisfying \eqref{condu} are gathered in the following table, where $0$ and $1$ stand for the convenient permutation or partitions on $[2]$:

\begin{center}
\begin{tabular}{|c|c|c|c|} \hline
$\alpha$ & $\beta$ &$\beta \alpha^{-1}$&$ A$\\\hline
0 &0&0& 1\\\hline
1&0&1&1\\\hline
0&1&1&1\\\hline
1&1&0&0 or 1\\\hline
\end{tabular}
\end{center}

The relative cumulants are given by (see \eqref{Wg}, \eqref{relcum})
\be
C_{1, 1} = -\frac{1}{n(n^2 -1)} \ , \ C_{0, 1} = - \frac{1}{n^2} + \frac{1}{n^2-1}=  \frac{1}{n^2(n^2-1)} \ , \ C_{0, 0} = \frac{1}{n^2}
\ee
The different products of traces are quite obvious, so that
plugging into \eqref{main2u}, we get
\begin{eqnarray} 
\label{cov1}
\lefteqn{
\kappa_2(T_{p,q}, T_{p',q'}) =} \\
&&  \frac{(p\wedge p')(q' \wedge q')}{n^2 -1} - \frac{(p \wedge p') q q'}{n(n^2-1)} - \frac{p p' (q\wedge q')}{n(n^2-1)} + \frac{pp'q'q'}{n^2(n^2-1)}. \nonumber 
\end{eqnarray}
In the limit $p/n \rightarrow s, q/n \rightarrow t, p'/n \rightarrow s', q'/n \rightarrow t'$, we get
\begin{eqnarray}
 \lim_n \kappa_2(T_{p,q}, T_{p',q'})  
  &=& (s\wedge s')(t\wedge t') - (s\wedge s')t t' - s s' (t\wedge t') + s s't t' \nonumber\\
&=& (s \wedge s' -s s')(t\wedge t - tt'). \label{cov2}
\end{eqnarray}
\subsubsection{The fourth cumulant}
We now give an estimate for  $\K_4(T_{p,q})$.
From \eqref{main2u} with $r=4$, 
\ben
 \label{sumk4}
\K_4 = \sum_{\alpha, \beta \in {\mathcal S_4}} \sum_{A} C_{\beta \alpha^{-1}, A} \!\ \Tr_{\alpha}(\bar D) \Tr_{\beta^{-1}}(D)
 \een
where $A$ runs over  the partitions of $[4]$ satisfying condition (\ref{condu}), 
  and finally
\[D_i = I_p \ , \ D_{\bar i}= I_q  \ \ (i \leq 4)\,.\]
We have now
\[\Tr_{\beta^{-1}}(D) = p^{\#(\beta)} \ , \ \Tr_{\alpha}(\bar D) = q^{\#(\alpha)}\,.\]
In \cite[Cor. 2.9]{Co}, Collins proved that the order of $C_{\beta \alpha^{-1},  A}$ is at most $n^{-8- \#(\beta \alpha^{-1})+2\#(A)}$. 
Finally, 
\ben
\label{bigO}
C_{\beta \alpha^{-1},  A}\Tr_{\beta^{-1}}(D) \Tr_{\alpha}(\bar D) =  O\left(n^{-8- \#(\beta \alpha^{-1})+2\#(A)}p^{\#(\beta)} q^{\#(\alpha)}\right)\,. \een
From equation (20) in \cite{mingo2007second}, we see that
\[2\#(A) + \#(\alpha) + \#(\beta) -\#(\beta \alpha^{-1})\leq 6\]
and  the expression in (\ref{bigO}) is of order
\begin{eqnarray*}
n^{-8- \#(\pi)+2\#(A)}p^{\#(\beta)} q^{\#(\alpha)}  &\leq& p^{\#(\beta)} q^{\#(\alpha)} n^{-2-\#(\alpha) - \#(\beta)}\\
&\leq& (p/n)^{\#(\beta) -1} (q/n)^{\#(\alpha)-1} pq n^{-4}\\ 
&\leq&  p^2q^2 n^{-4}\,.
\end{eqnarray*}
We conclude that
\ben \label{gather}
\K_4 = O\left(p^2q^2 n^{-4}\right)\,. 
\een

\subsection{Proof of Theorem \ref{maintheo}}
\subsubsection{Tightness}
According to Bickel and Wichura \cite[Theorem 3]{bickel1971convergence}, since our processes are null on the axes, the tightness of the family of distributions of $W^{(n)}$ is in force as soon as the condition ${\mathcal C}(\beta, \gamma)$ with $\beta >1$ is satisfied (see (2), (3) in \cite{bickel1971convergence}):
\ben \label{tight0}
\E(\vert W^{(n)}(B)\vert^{\gamma_1}\vert W^{(n)}(C)\vert^{\gamma_2}) \leq (\mu(B))^{\beta_1} 
 (\mu(C))^{\beta_2} 
 \een
 where $\gamma= \gamma_1 + \gamma_2 >0$ and $\beta= \beta_1+ \beta_2 >1$, $B$ and $C$ are two adjacent blocks in $[0,1]^2$ and $W^{(n)}(B)$ denotes the increment of $W^{(n)}$ around $B$,  given by
 $$ W^{(n)}(B) = W_{s',t'}^{(n)} -W_{s',t}^{(n)}- W_{s,t'}^{(n)}+ W_{s,t}^{(n)}$$
 for $B = ]s, s'] \times ]t,t']$, $\mu$ is a  finite positive  measure on $[0;1]^2$ with continuous marginals. \\
 From Cauchy-Schwarz inequality, \eqref{tight0} is implied by
 \ben \label{CS}
 \E(\vert W^{(n)}(B)\vert^{2\gamma_1}) \leq (\mu(B))^{2\beta_1} . 
 \een
Moreover, it is enough to consider blocks whose corner points are in $T^n= \{\frac{p}{n}, 0 \leq p\leq n\}
\times \{\frac{q}{n}, 0 \leq q \leq n\}$ (see \cite{bickel1971convergence}, p. 1665.)

\noindent
Let $p\leq p'\leq n$ and $q\leq q'\leq n$ and $B = ]\frac{p}{n}, \frac{p'}{n}] \times ]\frac{q}{n}, \frac{q'}{n}]$
\be  W^{(n)}(B) := \Delta^{(n)}_{p,q}(p',q') &=& Y_{p',q'}^{(n)} -Y_{p',q}^{(n)}- Y_{p,q'}^{(n)}+ Y_{p,q}^{(n)}\\
&=& \sum_{p+1 \leq i\leq p'}\sum_{q+1 \leq i\leq q'} |U_{i,j}|^2 - \mathbb{E}(|U_{i,j}|^2).\ee
If we  show that 
 there exists a constant $C$,  such that for all $n$,
 \ben
 \mathbb E \left[\left(\Delta^{(n)}_{p,q}(p',q')\right)^4\right]\leq C \frac{(p'-p)^2(q'-q)^2}{n^4},
 \een 
 then \eqref{CS} is satisfied with $\gamma_1 = 2, \, \beta_1 = 1$ and $\mu$ is the Lebesgue measure.
Since $\Delta^{(n)}_{p,q}(p',q')$ has the same distribution as $Y^{(n)}_{p'-p, q'-q}$, it is enough to show
\ben
\label{tight}
 \mathbb E \left[\left(Y^{(n)}_{p,q}\right)^4\right]
 = O(p^2q^2n^{-4})\,.
\een
If $X$ is a real random variable, an elementary computation gives
\ben\label{vark4}\mathbb E (X -\mathbb E X)^4 = \K_4 + 3 \K_2^2\,,\een
where $\K_r$ is the $r$-th cumulant of $X$.
Taking $X = T_{p, q}= \Tr D_1UD_{\bar 1}U^\star$, we saw above in (\ref{var0}) that 
\ben\label{varbis}\K_2 = \Var T_{p,q} \leq 2\frac{pq}{n^2}\,.\een
Gathering (\ref{vark4}) , (\ref{varbis}) and (\ref{gather}) we get that (\ref{tight}) is checked, which proves the tightness. 
\subsubsection{Finite-dimensional laws}
Let $(a_i)_{i \leq k} \in \mathbb R$, $(s_i, t_i)_{i \leq k} \in [0,1]^2$. We must prove the convergence in distribution of $X^{(n)} := \sum_{i=1}^k a_i W^{(n)}_{s_i,t_i}$ to a Gaussian distribution.\\
Let us denote $p_i = \lfloor ns_i\rfloor$, $q_i = \lfloor nt_i\rfloor$. Then
\begin{eqnarray*}
X^{(n)} &=& \sum_{i=1}^k a_i Y^{(n)}_{p_i,q_i} 
= \sum_{i=1}^k a_i [\Tr(D_{i} U D_{\bar i}U^\star) - \mathbb E(\Tr(D_{i} U D_{\bar i}U^\star))]
\end{eqnarray*}
where $D_{i} = I_{p_i}$, $D_{\bar i} = I_{q_i}$. \\
$\{D_{i}, D_{\bar i}, i=1, \ldots k \}$ are commuting projectors with a limit distribution 
$\{q_{i}, q_{\bar i}, i=1, \ldots k \}$  on a probability space $({\mathcal A}, \phi_1)$ with $\phi_1(q_{i})= s_i$, $\phi_1(q_{\bar i})= t_i$ and $q_iq_j = q_i$ if $u_i \leq  u_j$ (and $= q_j$ otherwise) where $u_i = s_i$ for $i$ odd and $u_i = t_i$ for $i$ even. \\
Let $r \geq 3$, then
\begin{eqnarray*}
\kappa_r(X^{(n)}, \ldots ,X^{(n)}) &= &\sum_{i_1, \ldots, i_r = 1}^k a_{i_1} \ldots a_{i_r} \kappa_r( Y^{(n)}_{p_{i_1},q_{i_1}},\ldots, Y^{(n)}_{p_{i_r},q_{i_r}}) \\
&=& \sum_{i_1, \ldots, i_r = 1}^k a_{i_1} \ldots a_{i_r} \kappa_r( X_{i_1}, \ldots, X_{i_r}) 
\end{eqnarray*}
where $X_{i_p} = \Tr(D_{i_p} U D_{\bar i_p} U^\star)$.
From Proposition \ref{mainprop2}
\ben
\lim_{n \rightarrow \infty} \kappa_r( X_{i_1}, \ldots, X_{i_r}) = 0.
\een
Now, the second cumulant is given by
$$\kappa_2(X^{(n)}, X^{(n)}) = \sum_{i,j=1}^k a_i a_j \kappa_2(\Tr(D_{i}UD_{\bar i}U^\star), \Tr(D_{j}UD_{\bar j}U^\star)).$$
From \eqref{cov2}
\begin{eqnarray*}
\lefteqn{ \lim_n \kappa_2(\Tr(D_{i}UD_{\bar i}U^\star), \Tr(D_{j}UD_{\bar j}U^\star)) } \\
&=& (s_i \wedge s_j -s_i s_j)(t_i\wedge t_j - t_it_j).
\end{eqnarray*}
Thus, we get the convergence of $X^{(n)}$ to a centered Gaussian distribution with variance
$$\sum_{i,j=1}^k a_i a_j  (s_i \wedge s_j -s_i s_j)(t_i\wedge t_j - t_it_j).$$
It follows that the finite-dimensional laws of the process $W^{(n)}$ converge to the finite-dimensional laws of the tied-down Brownian bridge. \QED

\section{The orthogonal case}
\subsection{Combinatorics for the orthogonal group}
Let ${\mathbb O}(n)$ denote the orthogonal group of size $n$ endowed with the Haar probability measure. 
 The generic element will be denoted by $O$ and its $(i,j)$ coefficient by $O_{ij}$. In \cite{CoSn}, Collins and Sniady proved the following integration formula on ${\mathbb O}(n)$, see also \cite[Theorem 2.1]{CoMa} for the following formulation. Let $G^{{\mathbb O}(n)}$ be 
 the Gram matrix\footnote{See footnote 1.} 
$$ G^{{\mathbb O}(n)} = (G^{{\mathbb O}(n)}(p_1, p_2))_{p_1,p_2 \in {\mathcal M}_{2k}} := (n^{{\rm loop}(p_1, p_2)})_{p_1, p_2\in {\mathcal M}_{2k}}\,,$$
where ${\mathcal M}_{2k}$ is the set of pairings of $[2k]$ defined in Section 2.1.
Then, the orthogonal Weingarten matrix $\W^{{\mathbb O}(n)}$ is the pseudo inverse of $G^{{\mathbb O}(n)}$. \\
\begin{prop} \label{propWeing-orth}
For every choice of indices ${\mathbf i} = (i_1, \ldots, i_{k}, i_{\bar 1}, \ldots, i_{\bar k})$ and  
${\mathbf j} = (j_1, \ldots, j_{k}, j_{\bar 1}, \ldots, j_{\bar k})$,
\begin{equation} \label{Weing-orth}
\mathbb E \left(  O_{i_1j_1} \ldots O_{i_{k}j_{k}}O_{i_{\bar 1}j_{\bar 1}} \ldots O_{i_{\bar k}j_{\bar k}}\right) 
 = \sum_{p_1, p_2 \in {\mathcal M}_{2k}} \delta_{\mathbf i}^{p_1}  \delta_{\mathbf j}^{p_2} \W^{{\mathbb O}(n)}(p_1, p_2)
\end{equation}
where 
 $\delta_{\mathbf i}^{p_1}$ (resp. $\delta_{\mathbf j}^{p_2})$ is equal to $1$  or $0$ if $\mathbf i$ (resp. $\mathbf j$) is constant on each pair of $p_1$ (resp. $p_2$) or not.
\end{prop}
We now identify ${\mathcal M}_{2k}$ as the quotient set ${\mathcal S}_{2k}/H_k$  where $H_k$ is a subgroup of ${\mathcal S}_{2k}$ known as 
the hyperoctahedral group and defined as follows (see \cite{GLM} and \cite{CC}). 
 With each  $g \in {\mathcal S}_{2k}$, we associate the product of disjoint transpositions:
$$ \eta(g) = \prod_{i=1}^k (g(i)\  g(\bar{i}) )\,,$$
which can be identified as an element of ${\mathcal M}_{2k}$. 
 
We set $\gamma =  \prod_{i=1}^k (i \ \bar{i})$ and define $H_k = \{ g \in  {\mathcal S}_{2k}, \gamma g = g \gamma\}$, the centralizer of $\gamma$. We have the following equivalence
$$ \eta(g) = \eta(g') \Longleftrightarrow \exists h \in H_k, g= g'h\,,$$
implying ${\mathcal M}_{2k} \cong {\mathcal S}_{2k}/H_k$.
 
According to Proposition 3.3 in \cite{CoSn}, see also \cite{CoMa}, $\W^{{\mathbb O}(n)}(\eta(g_1),\eta(g_2))$ depends only on the conjugacy class of $\eta(g_1)\eta(g_2)$ and we can define the orthogonal Weingarten function of ${\mathcal S}_{2k}$, denoted by $\Wl^{{\mathbb O}(n)}$ (see \cite[p. 511] {CC}) by
$$W\Lambda^{{\mathbb O}(n)}(g) = \W^{{\mathbb O}(n)}(\eta(Id), \eta(g))$$
and we have
$$\W^{{\mathbb O}(n)}(\eta(g_1), \eta(g_2)) = \Wl^{{\mathbb O}(n)}(g_1^{-1}g_2).$$
In the sequel, we shall drop the superscript ${\mathbb O}(n)$ keeping in mind for the asymptotics that 
$W\Lambda$ depends on the size $n$. \\
It is clear from the definitions that $\Wl$ is invariant on the classes of the double coset space $H_k\backslash {\mathcal S}_{2k}/ H_k$. \\
Then, formula \eqref{Weing-orth}  can be written as
\begin{equation} \label{Weing-orth2}
 \mathbb E \left(  O_{i_1j_1} \ldots O_{i_{k}j_{k}}O_{i_{\bar 1}j_{\bar 1}} \ldots O_{i_{\bar k}j_{\bar k}}\right)
  =  \frac{1}{|H_k|^2} \sum_{g_1, g_2 \in {\mathcal S}_{2k}} \delta_{\mathbf i}^{\eta(g_1)}  \delta_{\mathbf j}^{\eta(g_2)} \Wl(g_1^{-1}g_2).
 \end{equation}

We now  describe the generators of ${\mathcal S}_{2k} /H_k$ following  the presentation in \cite{GLM}, (see also \cite{CC} with a slightly different definition for particular permutations).
For $\epsilon= (\epsilon_1, \ldots, \epsilon_k)  \in \{-1, 1\}^k$, we define $\tau_\epsilon$ by 
$$\tau_\epsilon = \prod_{i, \epsilon_i= -1} (i \bar{i}) \in H_k.$$
For $\pi \in {\mathcal S}_{k}$, we define $t_\pi 
\in {\mathcal S}_{2k}$  by
$$t_\pi(i) = i; \; t_\pi(\bar{i}) = \overline{\pi(i)}\,.$$
We shall now parametrize ${\mathcal S}_{2k} /H_k$ using special permutations.
\begin{defn}
A pair $(\epsilon, \pi) \in \{-1,1\}^k \times {\mathcal S}_{k}$ is particular if, for any cycle $c$  of $\pi$, we have $\epsilon_i =1$ where $i$ is the smallest element of $c$. We say that the corresponding permutation 
 $g_{\epsilon, \pi}:= \tau_\epsilon t_\pi$ in ${\mathcal S}_{2k}$ is particular.
\end{defn}
\begin{prop} (see Theorem 8 in \cite{GLM}) \label{propLetac}
The class ${\mathcal S}_{2k} /H_k$ containing $\tau_\epsilon t_\pi$ has exactely $2^{\#(\pi)}$ elements of the form $\tau_{\epsilon'} t_{\pi'}$. Every class of ${\mathcal S}_{2k} /H_k$ contains exactely one particular permutation of the form  $\tau_\epsilon t_\pi$. There are $\frac{(2k)!}{2^k k!}$ particular permutations $g_l =  \tau_{\epsilon_l} t_{\pi_l}$, $l \leq \frac{(2k)!}{2^k k!}$.
\end{prop} 
Let $\Sigma \in {\mathcal S}_{2k}$. From the above Proposition, there exists a particular pair $(\epsilon, \sigma)$ such that $\Sigma = \tau_\epsilon t_\sigma h$ with $h \in H_k$. In particular,
$$ \Sigma \sim t_\sigma \, \mbox {in } \, H_k\backslash {\mathcal S}_{2k}/ H_k \; \rm{and} \; W\Lambda(\Sigma) = W\Lambda(t_\sigma).$$
We now recall how to find $\sigma$ from $\Sigma$ (see \cite[Proposition 18]{GLM}). First consider the pairing $\eta(\Sigma)$ and the 2-regular graph $\Gamma(\Sigma)$ (i.e. with all the cycles having even length) defined by 
$$\Gamma(\Sigma) = \eta(\Sigma) \cup \eta(Id).$$
Then, any cycle $\Gamma_j$ of length $2q_j$ of $\Gamma(\Sigma)$ is of the form 
$$ \Gamma_j = (i^{\pm}_{j1}, i^{\mp}_{j1}, \ldots, i^{\pm}_{jq_j}, i^{\mp}_{jq_j})$$
where the ordered couple $i^{\pm} i^{\mp}$ is such that $i^{\pm} i^{\mp} \in \{i, \bar{i} \}$ with $i^{\pm}= i$ implies $i^{\mp}= \bar{i}$ and $i^{\pm}= \bar{i}$ implies $i^{\mp}= i$. By convention, the starting point of $\Gamma_j $ is $i^{\pm}_{j1} = \min\{ \bar{i}_{jl}, l \leq q_j \}$. Set $\sigma_j = (i_{j1}, \ldots, i_{jq_j})$.
Then $\sigma = \prod \sigma_j$. 
\subsection{The variance of $T_{p,q}$}
From Proposition \ref{propWeing-orth} and the value of the Weingarten function for $k=2$ (\cite{CoSn}, \cite[Examples 2.1-3.1]{CoMa}, see also \cite{Ji}), we have: \\
For $i\not= k$, $j\not= \ell$,
$$ \mathbb E \left(O_{ij}^2 O_{k\ell}^2\right) = \frac{n+1}{n(n+2)(n-1)}. $$
For $j\not= k$
$$
\mathbb E \left(O_{ij}^2O_{ik}^2\right) = \frac{1}{n(n+2)}$$
and 
$$
\mathbb E O_{ij}^{4} = \frac{3}{n(n+2)}.
$$
From these 3 relations, we easily get
\begin{prop} \label{propvariance-orth}
The mean and the variance of $T_{p,q}$ are given by:
\begin{equation}
\label{mean-orth}
\mathbb E T_{p,q} = \sum_{i\leq p, j \leq q} \mathbb E O_{ij}^2 = pq\!\ \mathbb EO_{11}^2 = \frac{pq}{n}\,. 
\end {equation}
and 
\ben \label{var-orth}
\Var T_{p,q} =  2 pq \frac{n^2 - n(p+q) + pq}{n^2(n+2)(n-1)}. \een
Assume that $p/n \rightarrow s$, $q/n \rightarrow t$, then,
$$ \lim_n \frac{1}{n} \mathbb E T_{p,q} = st \ , \ \lim_n \Var T_{p,q}  =  2 st(1-(s+t) +st) = 2 st(1-s)(1-t).$$
\end{prop}

\subsection{Mixed moments of the variables $T_{p_i,q_i}$.} 
As in the unitary case, we need to compute cumulants $\kappa_r(X_1, \ldots X_r)$ where the variables $X_i$ are of the form
$$X_i = \Tr(D_iOD_{\bar{i}}O^{-1}).$$
We shall first establish an analogue of Proposition \ref{mainprop2}. We assume in the following that the matrices $D_i$ are deterministic and symmetric. Even if we shall deal later only with diagonal matrices, the computations are the same for general matrices.
\begin{prop} \label{mainprop2O}
Let $O$ be Haar distributed on ${\mathbb O}(n)$.
Let  $D= (D_1, \ldots D_k)$ and $\bar{D}= (D_{\bar{1}}, \ldots 
D_{\bar{k}})$ be two families  of deterministic and symmetric matrices of size $n$.
We set, for $ 1 \leq i \leq r$,
$$X_i = \Tr(D_iOD_{\bar{i}}O^{-1}).$$
Then,
\begin{equation} \label{main2}
\kappa_r(X_1, \ldots, X_r)= 
\sum_{ (\alpha, \beta, \epsilon) \in {\mathcal S}_{r} \times  {\mathcal S}_{r} \times \{\pm 1\}^r  } \lambda_{\alpha, \beta, r} \sum_{A}  C_{\sigma, A} \Tr_{\alpha}(D) \Tr_{\beta^{-1}}( \bar{D})
\end{equation}
where in the second sum
\begin{itemize}
\item
$\sigma \in {\mathcal S}_r$   is a function of $\alpha, \beta, \epsilon$ satisfying $t_{\alpha^{-1}} \tau_\epsilon t_{\beta} \sim t_\sigma \; {\rm in} \; H_k\backslash S_{2k} / H_k$ see \eqref{defsigma},
\item  $A \in {\mathcal P}(r)$  is such that $\sigma \leq A$ and 
  $A\vee \alpha \vee \beta= 1_r$,
\item   $C_{\sigma, A}$ are the relative cumulants of the orthogonal Weingarten function (see \cite{CoSn})
\item the combinatorial coefficient $\lambda_{\alpha, \beta, \epsilon}$ is 
$2^{r-\#(\alpha) - \#(\beta)}$.
\end{itemize}
\end{prop}
{\bf Proof:} 
We first give a formula for  the mixed moments, following \cite[Eq (20)]{CC}:
\begin{equation} \label{mixmom}
\E[\prod_{i=1}^k \Tr(D_iOD_{\bar{i}}O^{-1})] = \frac{1}{|H_k|^2} \sum_{g,g' \in {\mathcal S}_{2k}} W\Lambda (g^{-1} g') M_D^+(g) M_{\bar{D}}^- ((g')^{-1})
\end{equation}
where the $\pm$ generalized moments $M$ are functions on ${\mathcal S}_{2k}$, resp. $H_k$-right and $H_k$-left invariant: for $Y$ a $k$-tuple of random real matrices,
$$M_Y^+(g) = \E[\Tr_{\pi_l} (Y_1^{\epsilon_1(l)}, \ldots, Y_k^{\epsilon_k(l)})] \mbox{ if } g \in g_l H_k$$
$$M_Y^-(g) = \E[\Tr_{\pi_l} (Y_1^{\epsilon_1(l)}, \ldots, Y_k^{\epsilon_k(l)})] \mbox{ if } g \in  H_k g_l
$$
with $Y^{-1} = Y^t$ the transpose of $Y$. \\
Formula \eqref{mixmom} follows from \eqref{Weing-orth2} and straightforward computations.
Since the matrices $D$ are symmetric, the moments do not depend of the sequence $\epsilon$.
Using the parametrisation $g = \tau_{\epsilon_1} t_\alpha h$ and $g'= \tau_{\epsilon_2} t_\beta h'$ with $h, h' \in H_k$,
we can rewrite 
 \eqref{mixmom} as 
\begin{equation} \label{mixmom1}
\E[\prod_{i=1}^k \Tr(D_iOD_{\bar{i}}O^{-1})] = \sum_{\alpha, \beta \in {\mathcal S}_{k}, \epsilon \in\{\pm 1\}^k} 
\lambda_{\alpha, \beta,k} W\Lambda (t_{\alpha^{-1}} \tau_\epsilon t_{\beta}) \Tr_\alpha(D) \Tr_{\beta^{-1}}( \bar{D})
\end{equation}
The coefficient $\lambda_{\alpha, \beta,k}$ comes from the fact that we do not impose the sequences $\epsilon_1$ and $\epsilon_2$ associated with $\alpha$, resp. $\beta$ to be particular and $\epsilon = \epsilon_1 \epsilon_2$. From Proposition \ref{propLetac}, $\lambda_{\alpha, \beta, k} = 
2^{k-\#(\alpha)-\#(\beta)}$\\
As recalled before, with the triplet $(\alpha, \beta, \epsilon)$, we can associate a  particular pair $(\sigma, \epsilon')$ such that:
\begin{equation} \label{defsigma}
 t_{\alpha^{-1}} \tau_\epsilon t_{\beta} = \tau_{\epsilon'} t_\sigma h
\end{equation}
with $h \in H$. Then, $\Wl (t_{\alpha^{-1}} \tau_\epsilon t_{\beta}) = \Wl (t_{\sigma})$, denoted below by $\Wl (\sigma)$.
In the following, we will denote by $\sigma:= \sigma (\alpha, \beta, \epsilon)$ the permutation constructed above.

We now compute  the cumulant, following the same scheme as in \cite{mingo2007second}.
Let $C= \{V_1, \ldots V_k\}\in {\mathcal P}(r)$. We denote by $S_{V_i}$ the permutations on $V_i$.
\begin{eqnarray*}
\lefteqn{ \E_C(X_1, \ldots , X_r) = \prod_{i=1}^k \E \left(\prod_{j\in V_i} X_j\right)} \\
&&= \sum_{(\alpha_i, \beta_i, \epsilon_i) \in S_{V_i} \times S_{V_i} \times \{\pm 1\}^{|V_i|}; i \leq k} ( \prod_{i=1}^k \lambda_{\alpha_i, \beta_i, |V_i|}) W\Lambda(\sigma_1) \ldots W\Lambda(\sigma_k)\ldots \\
&& \qquad \qquad \qquad \qquad \dots \Tr_{\alpha_1}(D) \Tr_{\beta_1^{-1}}( \bar{D}) \ldots \Tr_{\alpha_k}(D) \Tr_{\beta_k^{-1}}( \bar{D}) \\
&&=\sum_{\begin{array}{c} (\alpha, \beta, \epsilon) \in {\mathcal S}_{r} \times  {\mathcal S}_{r} \times \{\pm 1\}^r \\
\alpha, \beta \leq C
  \end{array}} \lambda_{\alpha, \beta, r}
\Wl_C(\sigma) \Tr_{\alpha}(D) \Tr_{\beta^{-1}}( \bar{D})
\end{eqnarray*}
where $\sigma 
\in {\mathcal S}_r$ is 
  such that $\sigma \leq C$ 
 and (see \cite[section 3.4]{CoSn})
$$ \Wl_C(\sigma) = \prod_{i=1}^k \Wl(\sigma_{|V_i}).$$
Collins and Sniady defined the cumulants of orthogonal Weingarten functions that we will denote by $C_{\sigma, A}$  for $A$ a partition $\sigma$-invariant. They satisfy like (\ref{relcumU}) and (\ref{relcumrevU}):
$$\Wl_C(\sigma) = \sum_{A \in [\sigma, C]} C_{\sigma, A}\,.$$
Thus,
$$\E_C(X_1, \ldots X_r) = \sum_{\begin{array}{c} 
\alpha, \beta \leq C\\
 \epsilon \in  \{\pm 1\}^r
 \end{array}}
 \lambda_{\alpha, \beta, r}\sum_{A \in [\sigma, C]} C_{\sigma, A}
 \Tr_{\alpha}(D) \Tr_{\beta^{-1}}( \bar{D}).$$
 Now,
 
\begin{eqnarray*}
 \lefteqn{\kappa_r(X_1, \ldots, X_r) } \\
 &&= \sum_{C} \M (C, 1_r) \E_C(X_1, \ldots X_r) \\
 && = \sum_{C} \M (C, 1_r)  \sum_{\begin{array}{c} (\alpha, \beta, \epsilon) \in {\mathcal S}_{r} \times  {\mathcal S}_{r} \times \{\pm 1\}^r \\
 \alpha, \beta \leq C
   \end{array}} \lambda_{\alpha, \beta, r}\sum_{A \in [\sigma, C]} C_{\sigma, A}
 \Tr_{\alpha}(D) \Tr_{\beta^{-1}}( \bar{D})\\
 && = \sum_{ (\alpha, \beta, \epsilon) \in {\mathcal S}_{r} \times  {\mathcal S}_{r} \times \{\pm 1\}^r  }\lambda_{\alpha, \beta, r}\sum_{\sigma \leq A}  \sum_{\{C; A, \alpha, \beta  \leq C\}} \M (C, 1_r) C_{\sigma, A}
 \Tr_{\alpha}(D) \Tr_{\beta^{-1}}( \bar{D}) \\
 &&=  \sum_{ (\alpha, \beta, \epsilon) \in {\mathcal S}_{r} \times  {\mathcal S}_{r} \times \{\pm 1\}^r  } \lambda_{\alpha, \beta, r}\sum_{\sigma \leq A; A \vee \alpha \vee \beta = 1_r}  C_{\sigma, A} \Tr_{\alpha}(D) \Tr_{\beta^{-1}}( \bar{D})
 \ee
where 
 the last equality follows from
$$\sum_{C; A,\alpha,\beta \leq C} \M(C, 1_r) = 
\begin{cases}
&1 \;\;\;\;\mbox{if}\;   A \vee \alpha  \vee \beta = 1_r,\\
&0 \;\;\;\;\mbox{otherwise}.
\end{cases}
\qquad \Box$$
\medskip
 
We now study the asymptotic of \eqref{main2} when $n \rightarrow \infty$. We assume that the family of $n\times n$ matrices $(D_i)$  has a limit distribution.
It is known (see \cite[Theorem 3.16]{CoSn}) that the order of $C_{\sigma, A}$ is $n^{-2r- \#(\sigma) + 2\#(A)}
$. In \cite{CoSn}, the asymptotic of the cumulant is given in terms of a metric $l$ on pairings. We use that
$l(p_\sigma, p_{Id}) = |\sigma| := r - \#(\sigma)$ where $p_\sigma= \prod(i, \overline{\sigma(i)})$ is the pairing associated with $\sigma$ and $p_{Id}= \gamma$.

\noindent
Now, 
 $$  \Tr_{\alpha}(D) \Tr_{\beta^{-1}}( \bar{D}) = n^{\#(\alpha) + \#(\beta)}  \tr_{\alpha}(D) \tr_{\beta^{-1}}( \bar{D}) =  O(n^{\#(\alpha) + \#(\beta)}).$$
Therefore, for given $\alpha, \beta, \epsilon$, the corresponding term in the cumulant is of order
$$n^{-2r- \#(\sigma) + 2\#(A)+ \#(\alpha) + \#(\beta)},$$
where $\sigma$ and $A$ satify the conditions quoted in Proposition \ref{mainprop2O}, i.e.
\begin{equation} \label{condition}
\left\{  \begin{array}{l}
\sigma \mbox{ is the permutation defined by } \eqref{defsigma},\\
 \sigma \leq A \ \mbox{and}\ 
  A\vee \alpha \vee \beta= 1_r. \end{array} \right.
  \end{equation}

\begin{prop} \label{propcum}
Under the conditions \eqref{condition}, for $r \geq 3$,
\begin{equation} \label{neg}
-2r- \#(\sigma) + 2\#(A)+ \#(\alpha) + \#(\beta)+1\leq 0\,.
\end{equation}
\end{prop}

\begin{cor} 
\label{kappa3o}
Let $\{D_i,   i \in [2r]\}_n$ be a sequence  of deterministic and symmetric matrices of size $n$ which has a limit distribution and $X_i = \Tr(D_iOD_{\bar{i}}O^{-1})$. For $ r \geq 3$,
\begin{equation} \label{cumasymp}
\lim_{n\rightarrow \infty} \kappa_r(X_1, \ldots,X_r)=  0
\end{equation}
\end{cor}

  We first recall the following lemma (see \cite{mingo2007second} and the proof therein for the second assertion below).
  \begin{lem}   \label{sup}
  For $A, B \in {\mathcal P}(k)$ we have 
  $$\#(A) + \#(B) \leq k + \#(A\vee B)\,.$$
  Moreover, if there exists a block $A_i$ of $A$ and $B_j$ of $B$ such that $\#(A_i \cap B_j) = l$, then,
  $$\#(A) + \#(B) \leq k-l+1 + \#(A\vee B)\,.$$
\end{lem}

{\bf Proof of Proposition \ref{propcum}}
From the above lemma and Condition \eqref{condition}, we have:
\begin{equation} \label{O1}
\#(\alpha) + \#(\beta) \leq r + \#(\alpha \vee \beta)\,,
\end{equation}
\begin{equation} \label{O2}
\#(A) +  \#(\alpha\vee \beta) \leq r +1\,,
\end{equation}
\begin{equation} \label{O3}
\#(A)  \leq \#(\sigma)\,.
\end{equation}

The proof relies on the following property:
\begin{lem} \label{ptfix}
If $\alpha$ or $\beta$ has a fixed point, then there is a strict inequality in \eqref{O2} or in \eqref{O3}
\end{lem}
{\bf Proof of Lemma \ref{ptfix}} We denote by $\Sigma$ the permutation of ${\mathcal S}_{2r}$ given by 
$t_{\alpha^{-1}} \tau_\epsilon t_{\beta}$. 

1) Assume that $i$ is a fixed point of $\alpha$ and $\beta$. Then $(\Sigma(i), \Sigma(\bar{i})) = (i, \bar{i})$ or $(\bar{i}, i)$ depending on the sign of $\epsilon(i)$ and therefore, $i$ is a fixed point of $\sigma$. In this case,
$\sigma \vee \alpha \vee \beta \not= 1_r$ and therefore $\#(A) < \#(\sigma)$.

2) Assume that $i$ is a fixed point of  $\beta$ and $\alpha(i) \not= i$.
Then $\{\Sigma(i), \Sigma(\bar{i})\} = \{i, \overline{\alpha^{-1}(i)}\}$ and therefore the elements $\{i, \alpha^{-1}(i)\}$ are in the same cycle of $\sigma$ (thus in the same block of $A$), thery are obviously in the same block of $\alpha$. From Lemma \ref{sup}, the inequality is strict in \eqref{O2}.

3) Assume that $i$ is a fixed point of  $\alpha$ and $\beta(i) \not= i$. \\
$ \Sigma( \overline{\beta^{-1}(i)}) = i$ or $\bar{i}$ according to the sign of $\epsilon_i$ and 
$$ \Sigma( \beta^{-1}(i)) =  \beta^{-1}(i) \; \rm{or} \; \overline{\alpha^{-1}\beta^{-1}(i)}. $$
Thus, we find two distinct elements $(i,  \beta^{-1}(i))$ (or $(i, \alpha^{-1}\beta^{-1}(i))$) belonging to a cycle of $\sigma$ (thus to the same block of $A$) and belonging to a block of $\alpha \vee \beta$ (which is 
$\beta^{-1}$-invariant and  $\alpha^{-1}\beta^{-1}$-invariant). Therefore, the inequality \eqref{O2} is strict. $\Box$ \\

Let $\alpha$ and $\beta$ two permutations in ${\mathcal S}_{r}$. If $\#(\alpha) + \#(\beta) \leq r-2 + \#(\alpha\vee \beta)$, then, \eqref{neg} is satisfied, using \eqref{O2} and \eqref{O3}. From \eqref{O1}, it remains to 
 study the two 
  cases:

1) $\#(\alpha) + \#(\beta) = r + \#(\alpha\vee \beta)$,

2) $\#(\alpha) + \#(\beta) =  r-1 + \#(\alpha\vee \beta)$.

\noindent
{\bf Case 1)} It is not difficult to see that in this case, there exist two different fixed points for $\alpha$ or $\beta$. For example, a unique fixed point for both $\alpha$ and $\beta$ would imply:
$$ \#(\alpha) + \#(\beta)  \leq 1 + \frac{r-1}{2} + 1 + \frac{r-1}{2} = r+1.$$
This is not possible since in this case, $\alpha\vee \beta$ has at least two blocks.
So, let $i \not=j$ the two fixed points. Several situations can occur:

a) $i, j$ are  fixed points of both $\alpha$ and $\beta$. Then, as seen above, $i,j$ are fixed points of $\sigma$. In this case, $\#(A) \leq \#(\sigma) -2$, leading to \eqref{neg}.

b) $i$ is a fixed point of $\alpha$ and $\beta$ and $j$ is a fixed point of one of them. From the previous Lemma,
$$ \#(A) \leq \#(\sigma)-1 \;  \hbox{and} \; \#(A) +  \#(\alpha\vee \beta) \leq r-1 +1=r$$
leading to \eqref{neg}.

c) $i,j$ are fixed points of one of the two permutations. For example, $i,j$ are fixed points of $\beta$. If $\#\{i,j, \alpha^{-1}(i), \alpha^{-1}(j)\} \geq 3$, then, from Lemmas \ref{sup} and \ref{ptfix}, it is not difficult to see that $\#(A) +  \#(\alpha\vee \beta) \leq r-2  +1$ (either a block of $\sigma$ and a block of $\alpha$ have 3 common points or two blocks of  $\sigma$ and two blocks of $\alpha$ have two common points, etc.), leading to \eqref{neg}. \\
If $\#\{i,j, \alpha^{-1}(i), \alpha^{-1}(j)\} = 2$, in this case, $\alpha$ and $\sigma$ have a 
2-cycle $(i,j)$. This yields 
$\#(A) +  \#(\alpha\vee \beta) \leq r-1  +1$ but we also have $\#(A) < \#(\sigma)$ since 
$\sigma \vee (\alpha\vee \beta) \not= 1_r$. The other cases are similar, leading to
$$- \#(\sigma) + 2\#(A) +  \#(\alpha\vee \beta) \leq r -2 +1= r-1$$
proving \eqref{neg}.

\vspace{.3cm}
\noindent
{\bf Case 2)}: $\#(\alpha) + \#(\beta) =  r-1 + \#(\alpha\vee \beta)$. If there exists a fixed point, then one of the inequality in \eqref{O2} or \eqref{O3} is strict and we are done. 
If there is no fixed point, this implies that $r$ is even, $\#(\alpha) = \#(\beta) = \frac{r}{2}$ and $\alpha \vee \beta = 1_r$. An equality in \eqref{O2} is true only for $A= 0_r$, the partition in singletons. This would imply that $\sigma = Id$. Let $i \leq r$, 
$$\Sigma(i) = i \,  \hbox{or} \ \overline{\alpha^{-1}(i)} \; \hbox{and}\;  \Sigma(\bar{i}) =\beta(i)\, \rm{or} \ \overline{\alpha^{-1}\beta(i)}.$$
Therefore, $\sigma = Id$ implies that $\alpha^{-1} \beta(i) = i$ or $\alpha^{-1}(i) = \beta(i)$ and $\alpha$ and $\beta$ have a common 2-cycle. This is not possible ($\alpha \vee \beta =1_r$) except if $r= 2$. $\Box$.

\subsection{The fourth cumulant}
We give an estimate for $\kappa_4(T_{p,q})$. From the previous section,
$$\kappa_4(T_{p,q}) = \sum_{ (\alpha, \beta, \epsilon) \in S_4 \times  S_4 \times \{\pm 1\}^4  }\lambda_{\alpha, \beta, 4}\sum_{A; A \vee \alpha \vee \beta = 1_4}  C_{\sigma, A} \Tr_{\alpha}(D) \Tr_{\beta^{-1}}( \bar{D})$$ where
$D= (I_{p}, I_{p},I_{p},I_{p})$ and $\bar{D} = (I_{q}, I_{q},I_{q},I_{q})$. From the asymptotic behavior of the cumulant, each term is of order
\begin{equation} \label{ordrecum4}
n^{-8- \#(\sigma) + 2\#(A)}p^{\#(\alpha)} q^{\#(\beta)}
\end{equation}
where $\sigma, A, \alpha, \beta$ satisfy condition \eqref{condition}. \\
First assume that $\#(\alpha) = \#(\beta) =1$. Then the order of \eqref{ordrecum4} is at most:
$$ \frac{pq}{n^4} \leq \frac{p^2q^2}{n^4}.$$
Now assume that $\alpha$ or $\beta$ has at least two blocks.
From Proposition \ref{propcum}, the order of \eqref{ordrecum4} is at most
$$ n^{-1 - \#(\alpha) -\#(\beta)}p^{\#(\alpha)} q^{\#(\beta)}.$$
It is easy to see that this is at most of order $ \frac{p^2q^2}{n^4}$. For example, $\alpha$ and $\beta$ has at least two blocks,
$$n^{-1 - \#(\alpha) -\#(\beta)}p^{\#(\alpha)} q^{\#(\beta)} \leq \frac{p^2q^2}{n^5} \left( \frac{p}{n}\right)^{\#(\alpha)-2} \left( \frac{q}{n}\right)^{\#(\beta)-2}. $$ 
In the remaining cases, we find a majoration by $\frac{pq^2}{n^4}$ or  $\frac{p^2q}{n^4}$.
Therefore, we obtain
\begin{equation}
\label{kappa4o} \kappa_4(T_{p,q}) = O\left(\frac{p^2q^2}{n^4}\right)\,.
\end{equation}

\subsection{Proof of Theorem \ref{maintheo}}
The proof is similar to the proof in the unitary case, using the asymptotic vanishing cumulants of order $r \geq 3$ (Corollary \ref{kappa3o}) and the estimate (\ref{kappa4o}) for the fourth cumulant, ensuring the tightness of the family of distributions.

\section{Complementary remarks}

1) 
Since the sup norm is continuous for the Skorokhod topology, Theorem  \ref{maintheo} implies that 
\[\sup_{s,t \in [0,1]} |W^{(n)}(s,t)| \rightarrow \sup_{s,t \in [0,1]} |W^{(\infty)}(s,t)|\]
 in distribution, which implies that
\begin{equation}
\label{cvproba}
\sup_{s,t \in [0,1]} |\frac{1}{n}T_{\lfloor ns\rfloor, \lfloor nt\rfloor} - st| \rightarrow 0\end{equation}
in probability.

2) The definition of our process focusses on the trace of a random matrix $H_{p,q}$. 
This trace is a linear statistic of its empirical spectral distribution, i.e.
\[T_{p,q} = \Tr H_{p,q} = p \int x d\mu (x)\,,\]
where
\[\mu = \frac{1}{p} \sum_{k=1}^p \delta_{\lambda_k}\,,\]
and the $\lambda_k$'s are the eigenvalues of $H_{p,q}$.
 If we are interested  only in marginals ($p = \lfloor ns\rfloor, q= \lfloor nt\rfloor$, with $s,t \in (0,1)$ fixed), 
we can look directly at the asymptotic behavior of $\mu$ when $n \rightarrow \infty$ 
and then deduce results from the continuity of the mapping $\mu \mapsto \int xd\mu(x)$ on ${\mathcal M}_1([0,1])$. 
It is known that, if $p\leq q$ and $p+q\leq n$, the random matrix $H_{p,q}$ belongs to the 
Jacobi unitary/orthogonal ensemble (\cite{Collins} Theorem 2.2,  \cite{AGZ} Prop. 4.1.4), which entails 
that the joint distribution of eigenvalues has a density proportional to
\begin{equation}
\label{5.2}
\prod_{k=1}^p \lambda_i^{a-1} (1-\lambda_i)^{b-1}\prod_{1 \leq i< j\leq p} |\lambda_i- \lambda_j|^\beta
\end{equation}
where $a =(q-p+1)\frac{\beta}{2}$ and $b= (n -p-q+1)\frac{\beta}{2}$. 
The sequence of empirical spectral distributions converges to the generalized Kesten-McKay distribution. When $s \leq \min (t, 1-t)$, this distribution has a 
 density which can be parametrized by $s, t$ or by the endpoints of its support $(u_- , u_+)$ with $0\leq u_- <u_+\leq 1$:
\begin{equation}
\pi_{u_-,u_+}(x) = C_{u_-,u_+}\frac{\sqrt{(x-u_-)(u_+ -x)}}{2\pi x (1-x)}
\end{equation}
where
 $$C_{u_-,u_+}^{-1}:=\frac{1}{2}\left[1-\sqrt{u_-u_+}-\sqrt{(1-u_-)(1-u_+)}\right].$$
The relation between $(s, t)$ and $u_\pm$ is
\[u_\pm = \left[\sqrt{s(1-t)}\pm\sqrt{(1-s)t}\right]^2\,.\]
By continuity, in all cases, 
  we  recover a weak form of (\ref{cvproba}), i.e.
\[\lim_n \frac{1}{n}T_{\lfloor ns\rfloor, \lfloor nt\rfloor} = s\int x \pi_{u_-,u_+}(x) dx =st\,,\]
in probability.

It could  also be possible to recover the limiting fluctuations  of the marginal distribution of $T_{\lfloor ns\rfloor, \lfloor nt\rfloor}$ with $s,t$ fixed, i.e.
\[T_{\lfloor ns\rfloor, \lfloor nt\rfloor} - \mathbb E T_{\lfloor ns\rfloor, \lfloor nt\rfloor} \convlaw {\mathcal N}(0, \frac{2}{\beta}s(1-s)t(1-t))\]
from the known results on the fluctuations of linear statistics of $\mu$.  Actually, the result of Johansson \cite{Jofluc} 
is not specific of the Jacobi ensemble, but uses a model of random matrices invariant by conjugation, with polynomial external field.
Here, the ensemble is invariant but the potential is logarithmic (see (\ref{5.2})). 
 
The result is  a Gaussian limit with the good variance.

At another level, in the same asymptotics as above,  Hiai and Petz \cite{hiai2} proved that the family of empirical spectral distributions  
 satisfies the Large Deviation Principle in ${\mathcal M}_1([0,1])$ with speed $\beta n^{2}/2$ and good rate function, which in the case
 $s<t< 1/2$ is 
 \begin{eqnarray*}{\mathcal I}(\nu) &=& -s^2 \int\int \log |x-y| d\nu(x)d\nu(y) \\&&- s\int \left((1-s-t) \log (1-x) + (t-s) \log x\right) d\nu(x)  + I_0(s, t)\,.\end{eqnarray*}
where $I_0(s,t)$ is some constant (the limiting free energy).
Appealing again to the continuity of the mean, we deduce from the contraction principle that
$n^{-1} T_{\lfloor ns\rfloor, \lfloor nt\rfloor}$ satisfies the LDP at scale $n^{-2}$ with good rate function
\[I(c) = \inf  \{{\mathcal I}(\nu) ; \nu \in {\mathcal M}_1([0,1]) , \int_0^1 x d\nu(x) =c\}\,.\]
3) In multivariate (real) analysis of variance, the random variable $T_{p,q}$ is known as the Bartlett-Nanda-Pillai statistics.  
The exact distribution of $T_{p,q}$ is known by its Laplace transform which is an hypergeometric function of matrix argument (\cite{Muir} p.479).  Various asymptotic studies have been performed, essentially  $p,q$ fixed, $n \rightarrow \infty$ (large sample framework), or high-dimensional framework  with $q$ fixed, $n,p \rightarrow \infty$ and $p/n \rightarrow s < 1$ (see for instance \cite{fuji}). 
The asymptotic 
 regime of the present paper $(p/n \rightarrow s, q/n \rightarrow t)$ is considered in Section 4.4 of the book \cite{Baibook} and a CLT for the statistic $T_{p,q}$ may be deduced from Theorem 2.2 of \cite{Baietal}.
\bigskip
\bigskip 

\noindent{\bf Acknowledgement}
This work is supported by  the ANR project Grandes Matrices
Al\'eatoires ANR-08-BLAN-0311-01. C. D-M. thanks J.A.~Mingo for the references \cite{Redelmeier1} and \cite{Redelmeier2}. 

\renewcommand{\refname}{References}

\end{document}